\documentclass[runningheads]{llncs}
\usepackage{amsmath,amssymb}
\usepackage{algorithm,algorithmic}
\usepackage{pifont}
\usepackage{cases}
\usepackage{graphicx}
\usepackage{stackengine}    

\newcommand{\eps}{\varepsilon}

\newcommand{\one}{\mathbf{1}}

\newcommand{\R}{\mathbb{R}}

\newcommand{\bW}{{\bf W}}

\newcommand{\bw}{{\bf w}}

\newcommand{\ds}{\displaystyle}
\newcommand{\norm}[1]{\left\| #1 \right\|}

\newcommand{\braces}[1]{\left\{ #1 \right\}}

\usepackage[colorinlistoftodos,bordercolor=orange,backgroundcolor=orange!20,linecolor=orange,textsize=scriptsize]{todonotes}

\def\R{\mathbb R}

\def\one{{\mathbf 1}}

\def\<#1,#2>{\langle #1,#2\rangle}

\DeclareMathOperator*{\argmin}{arg\,min}
\newtheorem{assumption}[theorem]{Assumption}

\begin{document}
\title{An acceleration of decentralized SGD under general assumptions with low stochastic noise \thanks{The work of E. Trimbach and A. Rogozin was supported by Andrei M. Raigorodskii Scholarship in Optimization. The research of A. Rogozin is supported by the Ministry of Science and Higher Education of the Russian Federation (Goszadaniye) №075-00337-20-03, project No. 0714-2020-0005. This work started during Summer school at Sirius Institute.}}

\titlerunning{Catalyst acceleration of decentralized SGD}

%
%
\author{Trimbach Ekaterina \and
Rogozin Alexander}
\authorrunning{E.Trimbach \and A.Rogozin}
%
\institute{Moscow Institute of Physics and Technology}
\maketitle              
\begin{abstract}
Distributed optimization methods are actively researched by optimization community. Due to applications in distributed machine learning, modern research directions include stochastic objectives, reducing communication frequency and time-varying communication network topology. Recently, an analysis unifying several centralized and decentralized approaches to stochastic distributed optimization was developed in Koloskova et al. (2020). In this work, we employ a Catalyst framework and accelerate the rates of Koloskova et al. (2020) in the case of low stochastic noise.
\end{abstract}

\section{Introduction}\label{sec:introduction}

In this paper we consider an optimization problem with sum-type functional
\begin{equation}\label{eq:min_f_sum}
    f^{\star}  =  \min_{x \in \mathbb{R}^{d}}\left[f(x) := \frac{1}{n} \sum_{i=1}^{n} f_{i}(x)\right].
\end{equation}
Each $f_i$ defined in stochastic form
$$
f_{i}(x):=\mathbb{E}_{\xi_{i} \sim \mathcal{D}_{i}} F_{i}\left(x, \xi_{i}\right), 
$$
where $\xi_i$ is a random variable with distribution $\mathcal{D}_i$. Random variables $\braces{\xi_i}_{i=1}^n$ are independent and do not depend on $x$, as well. We seek to solve the problem \eqref{eq:min_f_sum} in a decentralized distributed environment, where each of $n$ agents locally holds $f_i$ and has access to stochastic gradients $\nabla F_i(x, \xi)$. The agents are connected to each other via a communication network.

In this paper, we accelerate  decentralized stochastic gradient descent under low noise conditions using the  Catalyst shell. The Catalyst approach was originally developed in \cite{catalystFirstOrder} for deterministic problems and generalized to a stochastic setting in \cite{genericAcceleration}. The Catalyst envelope allows to accelerate the convergence rate of a deterministic or stochastic optimization method with for strongly convex problems. We apply the Catalyst acceleration algorithm described in the article \cite{genericAcceleration} while extending it to the decentralized case. As a result, we manage to achieve acceleration under the condition of low stochastic noise.

In article \cite{koloskova2020unified} it is proved that complexity of Decentralized Stochastic Gradient Descent (DSGD) for computing such $X$ that $\mu\mathbb{E}\left\|\bar{x}-x^{\star}\right\|_2^{2} \leq \varepsilon$ is
$$
\tilde{\mathcal{O}}\left(\frac{L \tau}{\mu p} \log \frac{1}{\varepsilon} + \frac{\bar{\sigma}^{2}}{\mu n \varepsilon}+\frac{\sqrt{L}(\bar{\zeta} \tau+\bar{\sigma} \sqrt{p \tau})}{\mu p \sqrt{\varepsilon}}\right).
$$
In this article we propose an accelerated version od DSGD that achieves accuracy $\eps$ after
$$
\tilde{O}\left(\frac{\tau \sqrt{L}}{ p \sqrt{ \mu}}\log \frac{1}{\varepsilon}  +\frac{\sqrt{L} \bar \sigma^2}{n \mu \sqrt{\mu } \varepsilon} + \frac{\sqrt{L}(\bar{\zeta} \tau+\bar{\sigma} \sqrt{p \tau})}{\mu p \sqrt{\varepsilon}}\right) 
$$
iterations. In the case of low stochastic noise, i.e. when $\bar \sigma$ is small and the first summand is dominant, the total DSGD complexity decreases .

This paper is organized as follows. The rest of Section \ref{sec:introduction} is devoted to related work, notation and assumptions. After that, in Section \ref{sec:dsgd_overview} we overview the DSGD algorithm of \cite{koloskova2020unified} and provide its Catalyst-accelerated version in Section \ref{sec:accelerated_dsgd}. Finally, we cover the proofs of main theorems in Section \ref{sec:proofs}.

\subsection{Related work}
    
Decentralized optimization methods are based on applying gradient updates and communication procedures. Schemes based on direct distribution of gradient descent \cite{yuan2016convergence} and sub-gradient descent \cite{Nedic2009} have simple implementation, but only converge to a neighbourhood of the solution in case of constant step-sizes. Algorithms which converge to an exact solution are well developed in the literature, as well. Exact methods include EXTRA \cite{shi2015extra}, DIGing \cite{Nedic2017achieving} and NIDS \cite{li2019decentralized}. In many analyses the performance of a decentralized method depends on function conditioning $\kappa$ and graph condition number $\chi$, (the term $\chi$ typically characterizes graph connectivity). In order to enhance the dependence on $\kappa$ and $\chi$, accelerated schemes are applied. Accelerated methods can use either direct Nesterov acceleration (DNGD \cite{Qu2017}, Mudag \cite{ye2020multi}, OPAPC \cite{kovalev2020optimal}, Accelerated Penalty Method \cite{dvinskikh2019decentralized,li2018sharp}) or employ a Catalyst framework \cite{li2020revisiting}. Accelerated methods SSDA and MSDA \cite{scaman2017optimal} use dual conjugates to local functions $f_i$ held by the nodes. Moreover, decentralized schemes such as DIGing \cite{Nedic2017achieving} and Push-Pull Gradient Method \cite{Pu2018} are capable of working on a time-varying network. For more references and performance comparison see Table 1 in \cite{xu2020distributed}.

The methods mentioned above are deterministic, i.e. using gradient of local functions $f_i$. This paper is devoted to stochastic decentralized methods, which are particularly interesting due to their applications in distributed machine learning and federated learning \cite{konecny2016federated,mcmahan2016federated,mcmahan2017communication}. Several distributed SGD schemes were proposed and analyzed in \cite{koloskova2019decentralized,assran2019stochastic,tang2018d,lian2017can}. A Local-SGD framework was studied in \cite{zinkevich2010parallelized,stich2018local}. Finally, a recent paper \cite{koloskova2020unified} proposed a unified analysis covering many cases and variants of decentralized SGD.

\subsection{Notation and Assumptions}

Throughout the paper, we use small letters for vectors and capital letters for matrices. We also denote $\one = (1\ldots 1)^\top\in\R^n$. For matrix $X\in\R^{n\times d}$, we denote its rows $X = (x_1\ldots x_n)^\top$, the average of its rows $\bar x = 1/n \sum_{i=1}^n x_i$ and let $\bar X = X \frac{\one \one^\top}{n} = (\bar x \ldots \bar x)^\top$.

We also introduce standard assumptions for convex optimization.
\begin{assumption}\label{assum:L_smoothness} ($L$-smoothness)
Each function $F_i(x, \xi_i)$ is differentiable for every $\xi_i \in supp(\mathcal{D}_i)$  and L-smooth. In other words, there exists a constant $L \geq 0$ such that for every $x_1, x_2 \in \mathbb{R}^d$ it holds
\begin{equation}
\|\nabla F_i(x_1, \xi_i)  - \nabla F_i(x_2, \xi_i) \|_2 \leq L \|x_1 - x_2 \|_2
\end{equation}
\end{assumption}

\begin{assumption}\label{assum:mu_convexity} ($\mu$-strong convexity)
Every function $f_i: \mathbb{R}^d \to \mathbb{R}$ is strongly convex with constant $\mu \geq 0$. In other words, for every $x_1, x_2 \in \mathbb{R}^d$ it holds
\begin{equation}
f_i(x_1) - f_i(x_2) + \frac{\mu}{2}\|x_1  - x_2 \|_2^2 \leq \langle \nabla f_i(x_1) ,  x_1 - x_2 \rangle
\end{equation}

\end{assumption}

\begin{assumption}\label{assum:bounded_noise} (Bounded noise)\\
Let $\ds x^\star = \argmin_{x\in\R^d} f(x)$. Define 
$$
\zeta_{i}^{2}:=\left\|\nabla f_{i}\left(x^{\star}\right)\right\|_{2}^{2}, \quad \sigma_{i}^{2}:=\mathbb{E}_{\xi_{i}}\left\|\nabla F_{i}\left(x^{\star}, \xi_{i}\right)-\nabla f_{i}\left(x^{\star}\right)\right\|_{2}^{2} 
$$
and also introduce
$$
\bar{\zeta}^{2}:=\frac{1}{n} \sum_{i=1}^{n} \zeta_{i}^{2}
, \quad  \bar{\sigma}^{2}:=\frac{1}{n} \sum_{i=1}^{n} \sigma_{i}^{2}.
$$
We assume that $\bar{\sigma}^{2}$ and $\bar{\zeta}^{2}$ are finite.
\end{assumption}

In the decentralized case, computational nodes communicate to each other via a network graph. The communication protocol can be written using a matrix $\bW$, which elements $[\bW]_{ij}$ characterize the communication weights between individual nodes. Moreover, the network  changes over time, and at each communication round a new instance of $\bW$ is sampled from a random distribution $\mathcal{W}$ (for a detailed discussion of various cases of distribution $\mathcal{W}$, see \cite{koloskova2020unified}). For further analysis, we impose an assumption onto the mixing matrix, which is similar to that of \cite{koloskova2020unified}.


\begin{assumption}\label{assum:consensus} (Expected Consensus Rate)\\
The entries $w_{i j}^{(t)}$ of matrix $W^{(t)}$ are positive if and only if nodes $i$ and $j$ are connected at time $t$. Moreover, there exist two constants $\tau\in\mathbb{Z},~ \tau \geq 1$  and $p \in(0,1]$ such that for all integers $\ell \in$ $\{0, \ldots, T / \tau\}$ and all matrices $X \in \mathbb{R}^{d \times n}$ 
$$
\mathbb{E}_{\bW}\left\|X \bW_{\ell, \tau}-\bar{X}\right\|_{2}^{2} \leq(1-p)\|X-\bar{X}\|_{2}^{2}
$$
where $\bW_{\ell, \tau}=\bW^{((\ell+1) \tau-1)} \ldots \bW^{(\ell \tau)}$ ,  $\bar{X}:=X \frac{\mathbf{1 1}^{\top}}{n}$ and $\mathbb{E}$ is taken over the distribution $\bW^{t} \sim \mathcal{W}^{t}$.
\end{assumption}

\section{Decentralized SGD Algorithm}\label{sec:dsgd_overview}

Decentralized stochastic gradient descent allows to search for the minimum of a sum-type function by performing calculations in parallel on each machine. Decentralization increases the fault tolerance of the algorithm and ensures data security. Special cases of DSGD are SGD (one machine is used) and Local SGD (a central system is allocated that has access to all others nodes). Below we recall the DSGD algorithm and its main convergence result from \cite{koloskova2020unified}
\begin{algorithm}[H]
\caption{Decentralized SGD}
    \begin{algorithmic}[1]\label{alg:dsgd}
        \REQUIRE{Initial guess $X_0$, functions $f_i$,\\
        initialize step-size $\left\{\eta_{t}\right\}_{t=0}^{T-1}, $ number of iterations $T,$ \\
         mixing matrix distributions $\mathcal{W}^{t}$ for $t \in[0, T]$ ,\\ 
        for each i-th node initialize $x_{i}^{0} \in \mathbb{R}^{d}$ from $X_0$}
        \FOR{$t = 0, \ldots, T$}
            \STATE{Sample $\bW^{t} \sim \mathcal{W}^{t}$}
            \STATE{ Parallel processes for every task for worker $i \in[n]$}
            \STATE{Sample $\xi_{i}^{t}$ and compute stochastic gradient $g_{i}^{t}:=\nabla F_{i}\left(x_{i}^{t}, \xi_{i}^{t}\right)$}
            \STATE{$x_{i}^{t+\frac{1}{2}}=x_{i}^{t}-\eta_{t} g_{i}^{t} \quad $ }
            \STATE{$x_{i}^{t+1}:=\sum_{j \in \mathcal{N}_{i}^{t}} \bw_{i j}^{t} x_{j}^{t+\frac{1}{2}} \quad $ }
        \ENDFOR
    \end{algorithmic}
\end{algorithm}


\begin{theorem}\label{th:dsgd_convergence}
(Theorem 2 in \cite{koloskova2020unified}). Let Assumptions \ref{assum:L_smoothness}, \ref{assum:mu_convexity}, \ref{assum:bounded_noise}, \ref{assum:consensus} hold. Then for any $\varepsilon > 0$ there exists a step-size $\eta_t$ (potentially depending on $\varepsilon$) such that after running Algorithm \ref{alg:dsgd} for $T$ iterations $\varepsilon$-accuracy is attained in the following sense:
\begin{equation}\label{DSGD_accuracy}
    \sum_{t=0}^{T} \frac{w_{t}}{W_{T}}\mathbb{E}\left( f(\bar{x}^{t}) -f^{\star}\right) + 
    \mu \mathbb{E}\left\|\bar{x}^{T+1}-x^{\star}\right\|_2^{2} \leq  \varepsilon.
\end{equation}
Where   $w_{t}=(1-\frac{\mu}{2} \eta_t)^{-(t+1)}$,  $W_{T}:=\sum_{t=0}^{T} w_{t}$, $\bar{x}^{t}:=\frac{1}{n} \sum_{i=1}^{n} x_{i}^{t}$.
The number of iterations $T$ is bounded as
\begin{equation}\label{DSGD_complexity}
\tilde{\mathcal{O}}\left(\frac{\bar{\sigma}^{2}}{\mu n \varepsilon}+\frac{\sqrt{L}(\bar{\zeta} \tau+\bar{\sigma} \sqrt{p \tau})}{\mu p \sqrt{\varepsilon}}+\frac{L \tau}{\mu p} \log \frac{1}{\varepsilon}\right),
\end{equation}
where the $\tilde O$ notation hides constants and polylogarithmic factors.
\end{theorem}

\section{Accelerated DSGD}\label{sec:accelerated_dsgd}

\subsection{Catalyst shell}\label{subsec:catalyst_shell}

\subsubsection{Overview of Catalyst framework}

The Catalyst shell have gained a lot of attention recently, mainly due to their wide range of applications.
The method allows  to speed up the algorithms by wrapping them in a shell, in which at each step it is necessary to minimize some surrogate function. It was originally described in the article \cite{catalystFirstOrder} and \cite{catalystFromTheoryToPractice} and subsequently a huge number of adaptations of this algorithm recently were proposed.

This paper uses the implementation of the Catalyst envelope proposed in the article \cite{genericAcceleration}. Suppose we have some algorithm $\mathcal{M}$, which is able  to solve the problem
\begin{align}
F^{\star} = \min_{x \in \mathbb{R}^{d}} F(x)
\end{align} 
with some accuracy $\varepsilon$, where $F$ is $\mu$-strongly convex function. Paper \cite{genericAcceleration} suggests to choose a surrogate function $h_k$ satisfying the following properties:
\begin{align*}
    \left(\mathcal{H}_{1}\right)\quad &h_{k} \text{ is } (\kappa+\mu)\text{-strongly convex}. \\
    \left(\mathcal{H}_{2}\right)\quad &\mathbb{E}[h_{k}(x)] \leq F(x)+\frac{\kappa}{2}\left\|x-y_{k-1}\right\|^{2} \text{ for } x=\alpha_{k-1} x^{\star}+\left(1-\alpha_{k-1}\right) x_{k-1}. \\
    \left(\mathcal{H}_{3}\right)\quad &\forall \varepsilon_{k} \geq 0,~ \mathcal{M} \text{ can provide a point } x_{k} : \mathbb{E}\left[h_{k}\left(x_{k}\right)-h_{k}^{\star}\right] \leq \varepsilon_{k}\\
\end{align*}

\noindent We recall the Catalyst acceleration algorithm from \cite{genericAcceleration}.
\begin{algorithm}[H]
\caption{Generic Acceleration Framework with Inexact Minimization of $h_{k}$}
\begin{algorithmic}[1]\label{alg:generic_catalyst}
\REQUIRE{Input: $x_{0}$ (initial guess); $\mathcal{M}$ (optimization method); $\mu$ (strong convexity constant); $\kappa$ (parameter for $h_{k}) ; K$ (number of iterations); $\left\{\varepsilon_{k}\right\}_{k=1}^\infty$ (sequence of approximation errors).\\
Define $y_{0}=x_{0};~ q = \frac{\mu}{\mu+\kappa};~ \alpha_{0}=\sqrt{q}$ if $\mu \neq 0$.} 
\FOR{$k = 1, \ldots, K$}
    \STATE{Choose a surrogate $h_{k}$ satisfying $\left(\mathcal{H}_{1}\right),\left(\mathcal{H}_{2}\right)$ and calculate $x_{k}$ satisfying $\left(\mathcal{H}_{3}\right)  \text{ for } \varepsilon_{k};$ }
    \STATE{Compute $\alpha_{k}$ in (0,1) by solving the equation $\alpha_{k}^{2}=\left(1-\alpha_{k}\right) \alpha_{k-1}^{2}+q \alpha_{k} .$}
    \STATE{Update the extrapolated sequence $y_{k}=x_{k}+\beta_{k}\left(x_{k}-x_{k-1}\right)$ \\
    with $\beta_{k}=\frac{\alpha_{k-1}\left(1-\alpha_{k-1}\right)}{\alpha_{k-1}^{2}+\alpha_{k}}$.}
    \ENDFOR
    \RETURN Output: $x_{k}$ (final estimate).
    \end{algorithmic}
\end{algorithm}

\noindent The authors of \cite{genericAcceleration} provide the following result for Algorithm \ref{alg:generic_catalyst}.
\begin{theorem}\label{th:catalyst_rate}
After running Algorithm \ref{alg:generic_catalyst} for $k$ iterations, the following inequality holds:
\begin{align*}
    \mathbb{E}\left[F\left(x_{k}\right)-F^{\star}\right] \leq\left(1-\frac{\sqrt{q}}{2}\right)^{k}\left(2\left(F\left(x_{0}\right)-F^{\star}\right) + 4 \sum_{j=1}^{k}\left(1-\frac{\sqrt{q}}{2}\right)^{-j}\left(\varepsilon_{j}+\frac{\varepsilon_{j}}{\sqrt{q}}\right)\right)
\end{align*}
\end{theorem}

\subsection{Catalyst application}

In this section, we combine Algorithms \ref{alg:dsgd} and \ref{alg:generic_catalyst} and present a Catalyst-accelerated version of DSGD.

Consider a sequence of matrices $\{Y^{k}\}_{k=0}^\infty$, $Y^{k}\in \mathbb{R}^{d \times n}$ and matrix $X\in \mathbb{R}^{d \times n}$. Let $X = (x_1\ldots x_n)$, $Y^{k} = (y_1^{k}\ldots y_n^{k})$ and define a sequence of functions
\begin{equation}\label{eq:tilde_H_k_def}
 H^k (X)= \frac{1}{n}\sum_{1}^{n} \left[ f_i(x_i) +\frac{\kappa}{2}   \left\|x_i-y^{k-1}_i\right\|_2^{2} \right] - \frac{\kappa}{2}\sigma^{k-1}_y, 
\end{equation}
where
$$
    \sigma_y^{k-1} = \frac{1}{n} \sum_{i=1}^n \left\| y_i^{k - 1} \right\|_2^2 - \frac{1}{n^2} \left\|\sum_{i=1}^n y^{k-1}_i \right\|_2^2.
$$
Note that substituting $X = \bar{X} = X \frac{\mathbf{1 1}^{\top}}{n}$ yields
$$
 H^k(\bar X) = \frac{1}{n} \sum_{i=1}^n f_i(\bar x) + \frac{\kappa}{2} \left\| \bar x - \bar y ^{k-1} \right\|_2^2 =\\
f(\bar x) + \frac{\kappa}{2} \left\| \bar x - \bar y ^{k-1} \right\|_2^2.
$$
Also define
\begin{equation*}
 H_k^{\star} = \min_{x \in \mathbb{R}^{d}} \left[ \frac{1}{n}\sum_{i=1}^{n} \left[ f_i(x) +\frac{\kappa}{2}   \left\|x-y^{k-1}_i\right\|_2^{2} \right] - \frac{\kappa}{2}\sigma^{k-1}_y\right] =
\min_{x \in \mathbb{R}^{d}} \left[ f(x) + \frac{\kappa}{2} \left\|  x - \bar y ^{k-1} \right\|_2^2 \right].
\end{equation*}
Introduce $ h_k(x) = f(x) + \frac{\kappa}{2} \norm{ x - \bar y^{k-1}}_2^2$ and note that $H^k (\bar X) = h_k(\bar x)$.

\begin{algorithm}[H]
    \caption{Catalyst-accelerated decentralized SGD}
    \begin{algorithmic}[1]\label{alg:catalyst_dsgd}
        \REQUIRE{Number of outer iterations $K$. \\
        For every $k = 1, \ldots, K$, a step-size sequence $\{\eta_{t}\}_{t = 0}^{T_k-1}$.\\
        Define $q = \frac{\mu}{\mu + \kappa},~ \rho = \sqrt{q} / 3,~ \alpha_0 = \sqrt{q}$. Choose $\left\{\varepsilon_{k}\right\}_{k=1}^\infty$ (approximation errors) \\
        For each other $i$-th node initialize $x_{i}^{0} = y_i^{0} \in \mathbb{R}^{d}$.}
        \FOR{$k = 1, \ldots, K$}
            \STATE{Compute required number of iterations $T_k$ and step-size sequence $\{ \eta_{t}\}_{t=0}^{T_k-1}$ such that DSGD provides $\eps_k$ accuracy in $T_k$ iterations.}
            \STATE{Run DSGD on functions $h_k(x)$ for $T_k$ iterations using $x_i^{k-1}$ at node $i$ as initial guesses and setting step-size sequence to $\{ \eta_{t}\}_{t=0}^{T_k-1}$.}
            \STATE{Compute $\alpha_k\in (0, 1)$ by solving the equation $\alpha_k^2 = (1 - \alpha_k)\alpha_{k-1}^2 + q\alpha_k$.}
            \STATE{For each node, update $y_i^k = x_i^k + \beta_k(x_i^k - x_i^{k-1})$ with $\beta_k = \frac{\alpha_{k-1}(1 - \alpha_{k-1})}{\alpha_{k-1}^2 + \alpha_k}$.}
        \ENDFOR
        \RETURN matrix $X = \left[x_1^K, x_2^K, \ldots, x_n^K \right]$
    \end{algorithmic}
\end{algorithm}

\subsection{Convergence of Algorithm \ref{alg:catalyst_dsgd}}

Algorithm \ref{alg:catalyst_dsgd} can be written in matrix form
\begin{algorithm}[H]\label{alg:catalyst_dsgd_matrix}
    \caption{DSGD acceleration in matrix form.}
    \begin{algorithmic}[1]
        \REQUIRE{Number of outer iterations $K$, and constants $q = \frac{\mu}{\mu + \kappa}, \alpha_0 = \sqrt{q}$, sequence of approximation errors $\{ \varepsilon_k \}_{k=1}^{K}$. For each $i$-th node initialize $y_i^0 = x_{i}^{0} \in \mathbb{R}^{d}$.}
        \FOR{$k = 1, \ldots, K$}
            \STATE{Compute required number of iterations $T_k$ and step size $\{ \eta_{t_k}\}_{t=0}^{T_k-1}$ such that DSGD can return  $X_k$ : $\mathbb{E}\left[ H_{k}(\bar X)- H_{k}^{\star}\right] \leq \varepsilon_k$ after $T_k$ iterations.} 
            \STATE{ Compute $X_k$ =  DSGD ($X_{k-1},  H_k, T_k, \{ \eta_{t_k}\}_{t=0}^{T_k}$ ).}
            \STATE{ Compute $\alpha_k$ in (0, 1) by solving the equation $\alpha _k^2 = (1 - \alpha_k)\alpha_{k-1}^2 + q\alpha_k$. }
            \STATE{Update $Y_k = X_k + \beta_k(X_k - X_{k-1})$ with $\beta_k = \frac{\alpha_{k-1}(1 - \alpha_{k-1})}{\alpha_{k-1}^2 + \alpha_k}$.}
        \ENDFOR
        \RETURN matrix $X_K$
    \end{algorithmic}
\end{algorithm}

\begin{theorem}\label{th:acc_dngd_convergence}
Under Assumptions \ref{assum:L_smoothness}, \ref{assum:mu_convexity}, \ref{assum:bounded_noise}, \ref{assum:consensus}, there exist a number of outer iterations $K$, inner iteration numbers $\{T_k \}_{k = 0}^{K-1}$ and step-size sequences $\big\{\{ \eta_{t}\}_{t = 0}^{T_k-1}\big\}_{k=0}^{K-1}$, such that Algorithm \ref{alg:catalyst_dsgd} attains $\eps$-accuracy in the following sense: it yields $X^{T} \in \mathbb{R}^{d \times n}$ such that $\mu \left\| \bar x^{T} - x^{\star}\right\|_2^2 \leq \varepsilon$, where $\bar{x}^{T}:=\frac{1}{n} \sum_{i=1}^{n} x_{i}^{T}$. The total number of iterations $T$ is bounded as
$$
\tilde{O}\left(\frac{\tau \sqrt{L}}{ p \sqrt{ \mu}}\log \frac{1}{\varepsilon}  +\frac{\sqrt{L} \bar \sigma^2}{n \mu \sqrt{\mu } \varepsilon} + \frac{\sqrt{L}(\bar{\zeta} \tau+\bar{\sigma} \sqrt{p \tau})}{\mu p \sqrt{\varepsilon}}\right), 
$$
where the $\tilde{\mathcal{O}}$ notation hides logarithmic factors not dependent on $\varepsilon$.
\end{theorem}
Under low noise conditions, when $\bar\sigma$ is small, the dominant contribution is made by the first term, which results in an accelerated rate in comparison with Algorithm \ref{alg:dsgd}.

\section{Proofs of Theorems}\label{sec:proofs}

\subsection{Proof of Theorem \ref{th:dsgd_convergence}}

\begin{proof}

Theorem \ref{th:dsgd_convergence} is initially presented in the article \cite{koloskova2020unified} and proved in Lemma 15 of the Appendix of the same paper. To estimate the accuracy of the acceleration, we need to know exactly which variables are hidden under $\tilde{\mathcal{O}}$, so the following is a proof of Lemma 15, but analysed it in $\mathcal{O}$ notation.
\\
The proof of Lemma 15 of \cite{koloskova2020unified} shows that
\begin{align*}\label{eq:bound_from_lemma_15_koloskova}
\frac{1}{2 W_{T}} \sum_{t=0}^{T} w_{t} e_{t}+a r_{T+1} \leq \frac{r_{0}}{\eta} \exp [-a \eta(T+1)]+c \eta+64 B A \frac{\tau}{p} \eta^{2}
\end{align*}
where $r_{t}=\mathbb{E}\left\|\bar{x}^{(t)} - x^{\star}\right\|^{2},~ e_{t} = f\left(\bar{x}^{(t)}\right)-f\left(x^{\star}\right),~ a=\frac{\mu}{2},~ b=1,~ c=\frac{\bar{\sigma}^{2}}{n}$, $B=3L,~ A=\bar{\sigma}^{2}+\frac{18 \tau}{p} \bar{\zeta}^{2},~ d=\frac{96 \sqrt{3} \tau L}{p}$.

After that, the lemma considers several different values of step-size $\eta$.
\subsubsection{1) $\eta=\frac{\ln \left(\max \left\{2, a^{2} r_{0} T^{2} / c\right\}\right)}{a T} \leq \frac{1}{d}$}.\\
\\
In the case of $\eta=\frac{\ln \left(\max \left\{2, a^{2} r_{0} T^{2} / c\right\}\right)}{a T}$, the right-hand side of \eqref{eq:bound_from_lemma_15_koloskova} writes as\\
\begin{align*}
    & \frac{r_{0}}{\eta} \exp [-a \eta(T+1)]+c \eta+64 B A \frac{\tau}{p} \eta^{2} \\
    &\quad\leq \frac{r_{0} a T}{\ln(2)} \exp [-\ln(\max \left\{2, a^{2} r_{0} T^{2} / c\right\})]+\frac{c }{a T}\ln \left(\max \left\{2, a^{2} r_{0} T^{2} / c\right\}\right) + 64 BA \frac{\tau}{p} \ln^2 \left(\max \left\{2, a^{2} r_{0} T^{2} / c\right\}\right)\\
    &\quad\leq \frac{c}{\ln(2)a T} +\frac{c }{a T}\ln \left(\max \left\{2, a^{2} r_{0} T^{2} / c\right\}\right) + \frac{64 BA}{a^2 T^2} \frac{\tau}{p} \ln^2 \left(\max \left\{2, a^{2} r_{0} T^{2} / c\right\}\right)
\end{align*}

therefore to achieve accuracy $ \frac{1}{2 W_{T}} \sum_{t=0}^{T} w_{t} e_{t}+a r_{T+1}  \leq \varepsilon$  the number of iterations of the algorithm will be equal to $$\mathcal{O}\left( \frac{c}{aT} + \frac{1}{a}\sqrt{ \frac{BA\tau}{p \varepsilon}}\right)$$ \\
\\

\subsubsection{1) $\eta=\frac{1}{d} \leq \frac{\ln \left(\max \left\{2, a^{2} r_{0} T^{2} / c\right\}\right)}{a T}$}.\\
\\
In case $\eta=\frac{1}{d} \leq \frac{\ln \left(\max \left\{2, a^{2} r_{0} T^{2} / c\right\}\right)}{a T}$, the right-hand side of \eqref{eq:bound_from_lemma_15_koloskova} is estimated as follows.
\begin{align*}
& \frac{r_{0}}{\eta} \exp [-a \eta(T+1)]+c \eta+64 B A \frac{\tau}{p} \eta^{2} \\
&\quad\leq r_{0} d \exp \left[-\frac{a(T+1)}{d}\right] + \frac{c }{a T}\ln \left(\max \left\{2, a^{2} r_{0} T^{2} / c\right\}\right) + 64 BA \frac{\tau}{p} \ln^2 \left(\max \left\{2, a^{2} r_{0} T^{2} / c\right\}\right) \\
\end{align*}
Then to achieve accuracy $ \frac{1}{2 W_{T}} \sum_{t=0}^{T} w_{t} e_{t}+a r_{T+1}  \leq \varepsilon$
the number of iterations of the algorithm will be equal to  $$\mathcal{O}\left(\frac{d}{a}\ln{\frac{r_0 d}{\varepsilon}} + \frac{c}{a T} + \frac{1}{a}\sqrt{ \frac{BA\tau}{p \varepsilon}} \right)$$ . \\
\\
Hence to achieve accuracy $\sum_{t=0}^{T} \frac{w_{t}}{W_{T}}\left(\mathbb{E} f\left(\bar{\mathbf{x}}^{(t)}\right)-f^{\star}\right)+ \mu \mathbb{E}\left\|\bar{\mathbf{x}}^{(T+1)}-\mathbf{x}^{\star}\right\|^{2} \leq \varepsilon$ it is required to perform the following number of iterations:
\begin{align}
\mathcal{O}\left(\frac{\bar{\sigma}^{2}}{\mu n \varepsilon}+\frac{\sqrt{L}(\bar{\zeta} \tau+\bar{\sigma} \sqrt{p \tau})}{\mu p \sqrt{\varepsilon}}+\frac{L \tau}{\mu p} \log \frac{r_0 \tau L}{\varepsilon p}\right)
\end{align}
\end{proof}

\subsection{Proof of Theorem \ref{th:acc_dngd_convergence}}
Let us introduce a new function $h(x):~ \R^d\to \R$:
\begin{equation}\label{eq:h_k_def}
h_k(x) = f(x) + \frac{\kappa}{2} \left\|  x - \bar y ^{k-1} \right\|_2^2
\end{equation}
and note that  $ H(\bar X) = h(\bar x)$.

\subsubsection{Number of outer iterations\\}

Note that $h_k$ introduced in \eqref{th:acc_dngd_convergence} satisfies properties $\left(\mathcal{H}_{1}\right)$ and $\left(\mathcal{H}_{2}\right)$ in \cite{genericAcceleration}. Indeed, $h_k$ is $(\mu + \kappa)$-strongly convex and $\mathbb{E}\left[h_{k}(x)\right] \leq f(x)+\frac{\kappa}{2}\left\|x-y_{k-1}\right\|_2^{2}$. Therefore, Algorithm \ref{alg:catalyst_dsgd} is similar to Algorithm \ref{alg:generic_catalyst} in the article \cite{genericAcceleration}.

We choose the accuracy at step $k$ as $\varepsilon_k = \mathcal{O}\left((1-\sqrt{q} / 3)^{k}\left(f\left(\bar x_{0}\right)-f^{\star}\right)\right)$. The properties of $h_k$ allow to get an estimate on the number of outer iterations:
\begin{align}\label{outer_iter}
K=\mathcal{O}\left(\frac{1}{\sqrt{q}} \log \left(\frac{f\left(\bar x_{0}\right)-f^{\star}}{q \varepsilon}\right)\right),
\end{align}
where total accuracy is $\varepsilon = \varepsilon_K / q$.

For a detailed proof, see Section B.3 of \cite{genericAcceleration}.

\subsubsection{Number of inner iterations\\}

When solving the inner problem it is necessary to find such  $X_k$ that \\
$h_k(\bar x_k) - h_k^{\star} \leq \varepsilon_k$  by starting the algorithm DSGD from the point $X_{k-1}$. For further analysis, we define $L_h = L + \kappa,~ \mu_h = \mu + \kappa$.

The DSGD algorithm applied to minimizing $H_k$ for $T$ iterations guarantees accuracy
$$
\begin{array}{l}
\sum_{t=0}^{T} \frac{w_{t}}{W_{T}}\mathbb{E}\left(  H_k(\bar{X}^{t}) -  H_k^{\star}\right) + 
\mu_H \mathbb{E}\left\|\bar{x}^{T+1}-x^{\star}\right\|_2^{2} \leq \varepsilon_k.
\end{array}
$$
Using that $ H_k(\bar X) = h_k(\bar x)$, we can rewrite the previous statement as
$$
\begin{array}{l}
 \sum_{t=0}^{T} \frac{w_{t}}{W_{T}}\mathbb{E}\left( h_k(\bar{x}^{t}) -h_k^{\star}\right)+ 
\mu_H\mathbb{E}\left\|\bar{x}^{T+1}-x^{\star}\right\|_2^{2} \leq \varepsilon_k.
\end{array}
$$
This means that $\mathbb{E} \left\|  \bar x^{T+1} - x^{\star} \right\|_2 \leq  \frac{\varepsilon_k}{\mu_H}  $ and
$$\mathbb{E} \left(  h_k(\bar x^{T+1}) - h_k^{\star} \right) \leq\mathbb{E} \frac{L_h}{2}\left\|  \bar x^{T+1} - x^{\star} \right\|_2  \leq  \frac{L_h \varepsilon_k}{2\mu_h}.
$$
We conclude that if DSGD achieves accuracy 
$
\begin{array}{l}
\mu_H\mathbb{E}\left\|\bar{x}^{T+1}-x^{\star}\right\|_2^{2} \leq \varepsilon_k \\
\end{array}
$
after 
$$
\mathcal{O}\left(\frac{\bar{\sigma}^{2}}{\mu_h n \varepsilon}+\frac{\sqrt{L_h}(\bar{\zeta} \tau+\bar{\sigma} \sqrt{p \tau})}{\mu_h p \sqrt{\varepsilon}}+\frac{L_h \tau}{\mu_h p} \log \frac{r_{k - 1} \tau L_h}{\varepsilon p}\right) 
$$
iterations, then DSGD achieves accuracy 
$\mathbb{E} \left(   h_k(\bar x^{T+1}) - h_k^{\star} \right) \leq \varepsilon_k$  after 
$$\mathcal{O}\left(\frac{L_h\bar{\sigma}^{2}}{\mu_h^2 n \varepsilon_k}+\frac{L_h(\bar{\zeta} \tau+\bar{\sigma} \sqrt{p \tau})}{\sqrt{\mu_h}\mu_h p \sqrt{\varepsilon_k}}+\frac{L_h \tau}{\mu_h p} \log \frac{r_{k-1}  \tau L^2_h}{\mu_h \varepsilon_k p} \right)
$$ iterations.
\\

According to Proposition 5 in article
\cite{genericAcceleration}, for Algorithm 3 it holds 
\begin{align*}
\mathbb{E}\left[h_{k}\left(x_{k-1}\right)-h_{k}^{\star}\right]=\mathcal{O}\left(\varepsilon_{k-1} / q^{2}\right).    
\end{align*}
As a result, the inner complexity for each $k$-th step will be 

\begin{align*}
    &\mathcal{O}\left(\frac{L_h\bar{\sigma}^{2}}{\mu_h^2 n \varepsilon_k}+\frac{L_h(\bar{\zeta} \tau+\bar{\sigma} \sqrt{p \tau})}{\sqrt{\mu_h}\mu_h p \sqrt{\varepsilon_k}}+\frac{L_h \tau}{\mu_h p} \log \frac{r_{k-1} \tau L^2_h}{\mu_h \varepsilon_k p} \right) = \\
    &\mathcal{O}\left(\frac{L_h\bar{\sigma}^{2}}{\mu_h^2 n \varepsilon_k}+\frac{L_h(\bar{\zeta} \tau+\bar{\sigma} \sqrt{p \tau})}{\sqrt{\mu_h}\mu_h p \sqrt{\varepsilon_k}}+\frac{L_h \tau}{\mu_h p} \log \frac{(h_{k}\left(x_{k-1}\right)-h_{k}^{\star} ) \tau L^2_h}{\mu^2_h \varepsilon_k p} \right) = \\
    &\mathcal{O}\left(\frac{L_h\bar{\sigma}^{2}}{\mu_h^2 n \varepsilon_k}+\frac{L_h(\bar{\zeta} \tau+\bar{\sigma} \sqrt{p \tau})}{\sqrt{\mu_h}\mu_h p \sqrt{\varepsilon_k}}+\frac{L_h \tau}{\mu_h p} \log \frac{\varepsilon_{k-1} \tau L^2_h}{\mu^2_h \varepsilon_k p q^{2}} \right) = \\
    &\mathcal{O}\left(\frac{L_h\bar{\sigma}^{2}}{\mu_h^2 n \varepsilon_k}+\frac{L_h(\bar{\zeta} \tau+\bar{\sigma} \sqrt{p \tau})}{\sqrt{\mu_h}\mu_h p \sqrt{\varepsilon_k}}+\frac{L_h \tau}{\mu_h p} \log \frac{ \tau L^2_h}{\mu^2_h p q^{2}} \right)
\end{align*}

\subsubsection{Total complexity \\}
Let us sum the number of internal iterations by the number of external iterations. This yields total complexity:
\begin{align*}
T = \mathcal{O}\left(\sum_{k=1}^{K} T_k \right) 
&= \mathcal{O} \left( \sum_{k=1}^{K} \left(\frac{L_h\bar{\sigma}^{2}}{\mu_h^2 n \varepsilon_k} + \frac{L_h(\bar{\zeta} \tau+\bar{\sigma} \sqrt{p \tau})}{\sqrt{\mu_h}\mu_h p \sqrt{\varepsilon_k}}+\frac{L_h \tau}{\mu_h p} \log \frac{ \tau L^2_h}{\mu^2_h p q^{2}}\right)  \right) \\ 
&= \mathcal{O} \left( \sum_{k=1}^{K} \frac{L_h\bar{\sigma}^{2}}{\mu_h^2 n \varepsilon_k}+\sum_{k=1}^{K}\frac{L_h(\bar{\zeta} \tau+\bar{\sigma} \sqrt{p \tau})}{\sqrt{\mu_h}\mu_h p \sqrt{\varepsilon_k}}+\sum_{k=1}^{K}\frac{L_h \tau}{\mu_h p} \log \frac{ \tau L^2_h}{\mu^2_h p q^{2}} \right) \\
&= \mathcal{O} \left( \sum_{k=1}^{K} \frac{L_h\bar{\sigma}^{2}}{\mu_h^2 n \varepsilon_k}+\sum_{k=1}^{K}\frac{L_h(\bar{\zeta} \tau+\bar{\sigma} \sqrt{p \tau})}{\sqrt{\mu_h}\mu_h p \sqrt{\varepsilon_k}}+K\frac{L_h \tau}{\mu_h p} \log \frac{ \tau L^2_h}{\mu^2_h p q^{2}} \right)  
\end{align*}
After that, we recall that $\varepsilon_k = \mathcal{O}\left((1-\sqrt{q} / 3)^{K}\left(f\left(\bar x_{0}\right)-f^{\star}\right)\right)$ and therefore $\sum_{k=1}^{K}\frac{1}{\varepsilon_k}$ and $\sum_{k = 1}^{K}\frac{1}{\sqrt{\varepsilon_k}}$ are geometric progressions. Moreover, note that $q \leq 1$ and hence $\mathcal{O}\left(1 - \sqrt{1 - \frac{\sqrt{q}}{3} } \right) \geq \mathcal{O}\left(\sqrt{q}\right)$.

\begin{align*}  
T &= \mathcal{O} \left(  \frac{L_h\bar{\sigma}^{2}}{\mu_h^2 n \sqrt{q} \varepsilon_K}+\frac{L_h(\bar{\zeta} \tau+\bar{\sigma} \sqrt{p \tau})}{\sqrt{\mu_h}\mu_h p \sqrt{q\varepsilon_K}}+\frac{L_h \tau}{\mu_h p}\frac{1}{\sqrt{q}} \log \frac{ \tau L^2_h}{\mu^2_h p q^{2}}\log \left(\frac{f\left(\bar x_{0}\right)-f^{\star}}{q \varepsilon}\right) \right) \\  
&= \mathcal{O} \left(  \frac{L_h\bar{\sigma}^{2}}{\mu_h^2 n q^{3/2} \varepsilon}+\frac{L_h(\bar{\zeta} \tau+\bar{\sigma} \sqrt{p \tau})}{\sqrt{\mu_h}\mu_h p q\sqrt{\varepsilon}}+\frac{L_h \tau}{\mu_h p}\frac{1}{\sqrt{q}}  \log \frac{ \tau L^2_h}{\mu^2_h p q^{2}}\log \left(\frac{f\left(\bar x_{0}\right)-f^{\star}}{q \varepsilon}\right) \right).
\end{align*}

Let us choose $\kappa = L - \mu$, then $\mu_h = L,~ L_h = 2L - \mu,~ q = \frac{\mu}{L}$ and total complexity for accuracy $\mathbb{E} \left(f(\bar x^{T}) - f^{\star} \right) \leq \varepsilon$ will be 
\begin{align*}  
T &= \mathcal{O} \left(  \frac{\sqrt{L}\bar{\sigma}^{2}}{\mu \sqrt{\mu} n  \varepsilon}+\frac{\sqrt{L}(\bar{\zeta} \tau+\bar{\sigma} \sqrt{p \tau})}{\mu p \sqrt{\varepsilon}}+\frac{\sqrt{L} \tau}{\sqrt{\mu} p} \log \frac{L ^2\tau }{p\mu^2}\log \left(\frac{L(f\left(\bar x_{0}\right)-f^{\star})}{\mu \varepsilon}\right) \right) \\
&= \tilde{\mathcal{O}} \left(  \frac{\sqrt{L}\bar{\sigma}^{2}}{\mu \sqrt{\mu} n  \varepsilon}+\frac{\sqrt{L}(\bar{\zeta} \tau+\bar{\sigma} \sqrt{p \tau})}{\mu p \sqrt{\varepsilon}} + \frac{\sqrt{L} \tau}{\sqrt{\mu} p} \log \left(\frac{1}{ \varepsilon}\right) \right) 
\end{align*}
where $ \tilde{\mathcal{O}}$ -notation hides constants and polylogarithmic factors.

\bibliographystyle{splncs04}
\bibliography{references}

\end{document}